\documentclass[review]{elsarticle}

\usepackage{color}
\usepackage[normalem]{ulem} 

\journal{Energy Conversion and Management}

\begin{document}

\begin{frontmatter}

\title{Model predictive control of indoor microclimate: existing building stock comfort improvement}

\author[Skoltech]{A. Ryzhov}
\author[Skoltech]{H. Ouerdane \corref{mycorrespondingauthor}}
\cortext[mycorrespondingauthor]{Corresponding author}
\ead[url]{h.ouerdane@skoltech.ru}
\author[Skoltech]{E. Gryazina}
\author[Skoltech]{A. Bischi}
\author[MIT]{K. Turitsyn}

\address[Skoltech]{Center for Energy Systems, Skolkovo Institute of Science and Technology, 3 Nobel Street, Skolkovo, Moscow Region 143026, Russia}
\address[MIT]{Department of Mechanical Engineering, Massachusetts Institute of Technology, Cambridge, MA USA}

\begin{abstract}
Home retrofitting provides a means to improve the basic energy and comfort characteristics of a building stock, which cannot be renewed because of prohibitive costs. We analyze how model predictive control (MPC) applied to indoor microclimate control can provide energy-efficient solutions to the problem of occupants' comfort in a variety of situations principally imposed by external weather and room occupancy. For this purpose we define an objective function for the energy consumption, and we consider two illustrative cases: one building designed and built in recent times with modern HVAC equipment, and one designed and built several decades ago with poor thermal characteristics and no dedicated ventilation system. Our model includes various physical effects such as air infiltration and indoor thermal ''inertia mass'' (inner walls, floor, ceiling, and furniture), and also accounts for the impact of human presence essentially as heat and CO$_2$ sources. The influence on the numerical results of forecast horizons and of uncertainties due to inaccuracies in the weather and room occupancy forecasts, are analyzed. As we solve non-convex optimization problems using a linear and a nonlinear optimizer, full MPC performance is compared to both linearized MPC and a standard on/off controller. The main advantage of MPC is its ability to provide satisfactory solutions for microclimate control at the least possible energy cost for both modern and old buildings. As old buildings are usually not properly ventilated, we see  in light of our simulation results, that a supply ventilation system installation provides a solution for significant air quality improvement. 
\end{abstract}
\end{frontmatter}

\section{Introduction}
Buildings designed and built during the 20th century and before form the largest part of those in use today. Their basic characteristics in terms of energy efficiency, environmental impact and comfort are de facto below par with the present standards driven by low-carbon policies, sustainability and features of modern energy systems. Furthermore, aging of buildings is accelerated by poor/harsh weather conditions \cite{Sadauskiene2009,Liisma2015}, favoring wall degradation because of, e.g., penetration of moisture in cracks and cavities of the walls and around windows, which in turn enhances heat leaks. In fact, buildings that host people in general, either public offices or dwellings, pose particular design and engineering challenges, as the human presence entails particular requirements in terms of comfort and use of energy. 

Generally speaking, comfort is usually characterized by a variety of ranges of measurable quantities such as temperature, humidity, concentration of CO$_2$ generated by occupants, other pollutants and particulate matter (PM), as well as lighting, and radiant heat. Given the complexity and interplay among these quantities, comfort constitutes an active field of research, which drives innovative design in the building sector by setting bounds on them \cite{Ortiz2017,ASHRAE}: indoor comfort temperatures typically lie between 19$^{\circ}$C and 25$^{\circ}$C, the relative humidity within the 45\% to 65\% range, and ideally CO$_2$ concentration should be kept to a minimum comparable to, e.g., the environmental level which is typically below 400 ppm \cite{MacFarlingMeure2006}. The standards set by the Russian government for microclimate in living apartments are more flexible \cite{RussianStandards}: for temperature, the range is from 18$^{\circ}$C to 26$^{\circ}$C, for humidity from 30\% to 65\%, and for ventilation there is no mention of a CO$_2$ concentration level, but rather the number of times air should be replaced in a room. Note that CO$_2$ concentration deserves particular attention: levels recorded in houses may reach values as high as 5000 ppm, which dramatically affects people's behavior as recently shown in dedicated studies \cite{Satish2012,Allen2016}; in fact, even at moderate levels of concentrations, say 2500 ppm, routinely reached in classes and conference rooms, intellectual activities like strategic thinking, initiative, information utilization, are negatively impacted. 

Energy requirements of occupied buildings amount to 40\% of primary energy consumption, 70\% of electricity consumption, and 30\% of greenhouse gas emission \cite{Levine2007} with most of the energy consumption being due to heating, ventilation and air conditioning (HVAC) systems operation \cite{Nejat2015}. So, as traditional microclimate control solutions, like opening a window, are energy-intensive, a proper design and control of HVAC systems is crucial both from the energy saving and health viewpoints. Further, while the design of new buildings should account for their shape and the surface area of their envelope to minimize, so that energy efficiency is enhanced \cite{AlAnzi2009,Hwang2011}, a proper use of various tools used for commissioning in buildings \cite{Ginestet2013} and consideration of renewable energy options as solutions \cite{Rezaie2013} can provide a better control of energy consumption. 

For the current large stock of existing buildings their shape is given and difficult to modify, so efficiency improvement can be reached only at moderate costs by the upgrading of the thermal insulation and the installation of recuperation systems for ventilation as well as the optimal control of microclimate. By doing so, the energy consumption can be minimized while maintaining an acceptable level of comfort \cite{Herrando2016,Zubiaga2016}. However, besides hardware (sensors, controllable actuators, information exchange between zones, computational facilities installed in every HVAC unit) challenges from the control viewpoint are many and include model accuracy, stochasticity in internal heat sources, and quality of forecast of external weather conditions. 

In the present work, we are interested in indoor comfort through an accurate control of the smart systems setting the indoor microclimate with a focus on both high- and low-efficiency buildings. To that end, we use a model predictive control approach \cite{MorariWeForc12}. We aim to see how a sizable decrease of energy consumption may be efficiently achieved under various conditions: building characteristics, occupancy, and weather, while the desired level of comfort is maintained. More precisely, we consider a state-of-the-art 100 m$^2$ conference room located in a modern office building near Moscow (test case 1) and the typical urban Russian panel building series P44 (test case 2). 

Most of Russian buildings are heated by district heating (DH) systems. If no access to DH exists, individual gas boilers are used. Some radiators are equipped with valves that can be manually tuned to control the heating level. Air conditioners are widely used in regions with warm summers and are controlled mainly by setting a desired temperature on on/off controllers. Natural ventilation dominates in dwellings which implies that windows must be opened even in winters despite its inefficiency and loss of comfort. Commercial buildings are equipped with air handling units (AHU) mostly based on constant air volume (CAV) technology resulting in uncontrollable ventilation rate while variable air volume (VAV) systems, which are more energy-efficient and provide a higher level of control of temperature and humidity, are not so much used because of their price.

Test cases 1 and 2, the detailed characteristics of which will be detailed in the main text, correspond to buildings submitted to the same weather conditions over the years, and harsh winters in particular; so their comparison provides quantitative information on the actual added value of the state-of-the-art technology and design, and also on the model predictive control approach we adopt. The article is organized as follows. In Section 2 we briefly review the main methods which can be applied to microclimate control, and discuss our choice for MPC in the context of the present work. In Section 3, we present our approach to the comfort and energy saving problem in buildings, based on the linear and model predictive control methods. In section 4, we implement numerically our approach to two illustrative test cases, and analyze and discuss the obtained results, considering also the possible presence of a supply ventilation system (SVS) in homes. The article ends with concluding remarks.

\section{\label{model}Overview of the main control methods} 
Different control methods have been developed for HVAC systems, including rule-based, model predictive control \cite{Killian2016}, PI/PID and artificial neural networks-based approaches~  \cite{Afram2014,Wang2017,Hazyuk2014}. In Russia, on/off controllers are usually used in dwellings, while rule-based ones are used in commercial buildings.

\subsection{Rule-based control} 
Rule-based control is widely used for temperature control on building scale and also taken as industry standard. Here all control inputs are taken from rules of the kind ``if \textit{condition} then \textit{action}'' \cite{report}. The main drawback of the method is that it needs tuning during operation with conditions changing to provide optimality. Unfortunately, this type of control may lead to temperature overshoots or synchronization effects \cite{GuptaCDC12} and hence may not provide optimal solutions from energy consumption viewpoint. Moreover, rule-based control strategies do not use predictions to perform better control actions.

\subsection{PI and PID controllers}
Various controllers for thermostats improve transient dynamics to reach the objectives set. Classical proportional-integral (PI), proportional-integrate-derivative (PID) controllers as well as fuzzy, adaptive, neural network  controllers \cite{Dounis2009} are capable to ensure that the desired temperature is reached and to exploit HVAC for providing demand response services \cite{HiskensDemandResponse16}. But the controllers performance is sensitive to the choice of the gains. Having tuned controller parameters for a particular objective like following a certain temperature, the designer looses the ability to reach other objectives such as minimization of an energy consumption, as well as to find a compromise among them. Besides, as the control goal is to maintain the prescribed air characteristics, comfort requirements dominate energy saving objectives. In other words PI and PID controllers cannot ensure optimal control or stability. 

\subsection{Model predictive control}
Model predictive control \cite{MPCbook} provides a suitable framework for microclimate control as it ensures optimality using an appropriate strategy. Instead of operating on an infinite time horizon as in classical optimal control, MPC yields optimal results for finite prediction horizons. At each time-step, MPC solves an optimal control problem over a given prediction horizon and obtains the control parameters and states that satisfy both the dynamics and constraints. Finally, it synthesizes a control signal that minimizes objective functions (operation cost, energy consumption) and satisfies comfort constraints. Existing methods for predicting the energy consumption of households, commercial, industrial and municipal consumers, using mostly regression models, do not take into account the new technological capabilities of consumers and do not allow to predict a demand response in the event of a significant change in pricing policy. From a control perspective, MPC is capable to face all the requirements and possesses additional useful features such as: ability to use occupancy profile and incorporate weather forecast, moderating trade-off between preferences that may vary between thermal comfort and energy efficiency covering all intermediate solutions.

For MPC the formulation of an accurate microclimate dynamics model is a crucial problem \cite{Privara1stEx,HazyukPart1,HazyukPart2}. Linear models allow to use linear optimization solvers which are fast enough, but they provide suboptimal solutions due to lower accuracy. Bilinear models can also capture the airflow effect of ventilation: the heat flux is proportional to the mass flow rate (control variable) and the temperature difference (state variable) thus improving the model accuracy. In turn, sequential linear programming \cite{Storen95} solves bilinear optimization problems fast enough. The MPC approach has been implemented successfully for several problems pertaining to thermal comfort and microclimate control \cite{MorariCDC10,Sturzenegge2015,CastillaThermal14,Scherer2014}.

\subsection{Centralized and distributed control schemes}
MPC can also be used for multiple zones control as both centralized and decentralized control schemes exist \cite{Scherer2014}. The centralized ones are difficult to implement in practice due to the significant optimization variables space needed, high number of communications and sensitivity to local changes of HVAC. To avoid these difficulties decentralized versions are designed and implemented. For example, the system can be decomposed into several control systems. Local subsystems may operate autonomously (purely decentralized scheme) or exchange information between neighboring zones (distributed scheme) and send limited information to the higher-level system \cite{Afram2017,Zhang2017}. 

\section{Formulation of the problem} 
\begin{figure}
	\centering
	\includegraphics[width=120 mm]{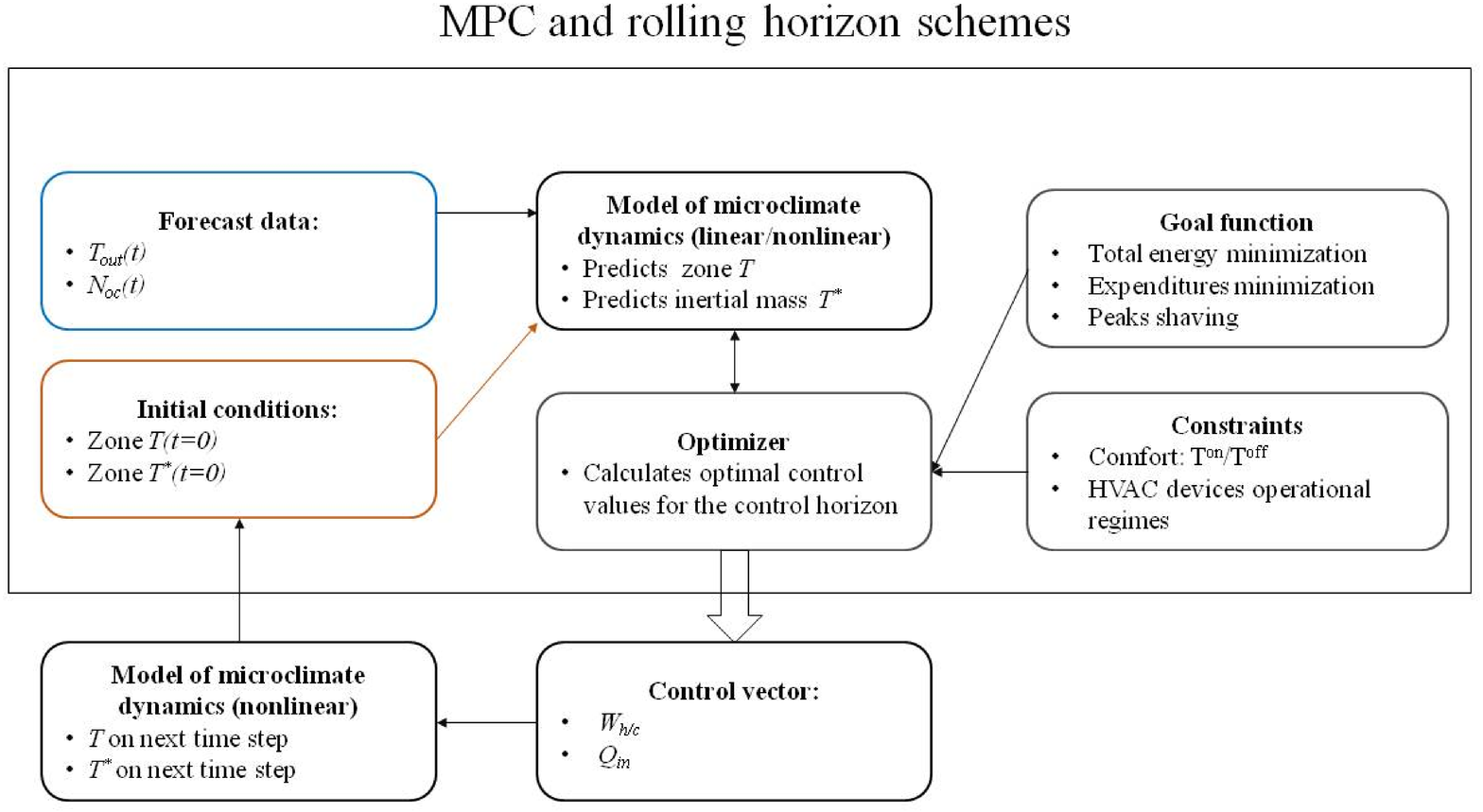}
	\caption{MPC and rolling horizon schemes.}
	\label{fig:Scheme}
\end{figure}

In the present work, we apply MPC to microclimate control in rooms equipped with DH radiators, air conditioners and an inflow ventilation. We simulate one generic day for three types of weather: cold, mild, and hot as defined further below. We thus set as an initial condition an ``effective temperature'', which amounts to neglecting the exact details of previous history that may have an impact if one wants to consider cases of sudden change of weather. In practice the effective temperature depends on the building usage/external conditions (absence days/ week ends / different weather etc.), thus leading to the simulation of different scenarios that may provide a set of initial conditions for a 1-day simulation and hence to a sensitivity analysis of the controller performance on the initial conditions, which is beyond the scope of the article. Note that the practical or experimental implementation of the controller does not allow to explicitly measure the effective temperature, but a measured room temperature gradient determined by thermal inertia, which can be evaluated during periods with the least possible amount of unpredicted factors, can instead be provided. This is very close to the approach we use, providing the model with an initial value of the effective temperature (21$^{\circ}$C in the present work.) The scheme of MPC and its testing using rolling horizon simulation is shown in Fig.~\ref{fig:Scheme} and explained as follows.

Our numerical simulations follow cycles, which first are based on the MPC block that consists of six logical parts necessary: the ``Forecast data'' part provides the outside temperature $T_{\rm out}$ and occupancy $N_{\rm oc}$ for the next 24 hours; the ``Initial conditions'' part provides the initial condition data to solve Eqs.~(\ref{eq:phys1}) and (\ref{eq:phys2}) below, for the microclimate dynamics: air temperature $T$, effective temperature of inertia mass $T^\star$, and CO$_2$ concentration taken to be equal to the environmental level in the beginning; note that when practically implementing MPC for HVAC control, the initial data correspond to the filtered temperature and CO$_2$ concentration values given by sensors installed in a zone; the ``Model of microclimate dynamics'' part provides the solutions ($T$, $T^\star$, CO$_2$) of Eqs.~(\ref{eq:phys1}) and (\ref{eq:phys2}), at the next time step based on forecast data, initial conditions, heating/cooling power $W_{\rm h/c}$, ventilation mass flow rate $Q_{\rm in}$; the ``Optimizer'' finds optimal values of $W_{\rm h/c}$ and $Q_{\rm in}$ accounting for both the ``goal function'' and ''constraints'' parts. The result of the MPC cycle is the ``Control vector'' containing the $W_{\rm h/c}$ and $Q_{\rm in}$ values at each control time step, taken as 1~hour in our cases. The control vector is then fed to the rolling horizon component of our simulation scheme, which undertakes a cycle external to the MPC block; when completed, we use the obtained microclimate data as updated initial conditions to run the next MPC cycle. The procedure is repeated until the simulated time reaches 24 hours. 

\subsection{Model of a room microclimate} 
The central part of the method is the microclimate physical model, which provides the time evolution of, e.g., the temperature and CO$_2$ concentration from the knowledge of the current state and action. In turn, this allows us to take into account the microclimate dynamics in optimal actions calculation which is the main advantage of MPC. From the principles of mass and energy conservation applied to a microclimate model \cite{Sturzenegge2015,Rincon2016}, we obtain the following set of differential equations:

\begin{eqnarray}
\nonumber
\label{eq:phys1}
mC_p \frac{{\rm d}T}{{\rm d}t} & = & U (T_{\rm out} - T) +U^{\star} (T^{\star} - T) +W_{\rm oc}N_{\rm oc}\\
&+&W_{\rm h/c} + C_p Q_{\rm in} (T_{\rm in} - T) + C_p m R_{\rm r} (T_{\rm out} - T)\\
m^{\star} C^{\star} \frac{{\rm d}T^{\star}}{{\rm d}t} & = & -U^{\star} (T^{\star} - T)
\label{eq:phys2}
\end{eqnarray}

\noindent where $t$ [h] is the time (which may be converted into seconds for some calculations); $W_{\rm h/c}(t)$ [W] is the time-dependent heating/cooling power; $Q_{\rm in}(t)$ [kg$\cdot$s$^{-1}$] is the time-dependent ventilation mass flow rate; $m$ [kg] is the mass of air in the room derived from the ideal gas equation of state: $m=PV/(RT)$, with $P=10^5$ Pa being the normal atmospheric pressure, $V$ the room volume, $R=287.03$ J$\cdot$kg$^{-1}\cdot$K$^{-1}$ is the value of the ratio of the ideal gas constant (8.314 J$\cdot$mol$^{-1}\cdot$K$^{-1}$) and the air molar mass (0.028966 kg$\cdot$mol$^{-1}$), and $T$ [$^{\circ}$C] the indoor air temperature (converted into kelvin for some calculations); $C_p=1$ kJ$\cdot$kg$^{-1}\cdot$K$^{-1}$ is the specific heat capacity of air at constant pressure; $U$ [W$\cdot$K$^{-1}$] is the linear heat transfer coefficient between the indoor and outside air at temperatures $T$ and $T_{\rm out}$ [$^{\circ}$C] respectively; $R_{\rm r}$ [h$^{-1}$] is the rate of air replaced by infiltrated air per one hour \cite{Solupe2014}; $m^{\star}$ [kg] is the ``inertia mass'', i.e. the accumulated mass of the walls, floor, ceiling and furniture of the room with $T^{\star}$ [$^{\circ}$C] being its effective temperature and $C^{\star}$ being its effective average heat capacity [J$\cdot$K$^{-1}$]; $U^{\star}$ [W$\cdot$K$^{-1}$] denotes the linear heat transfer coefficient between the inertia mass and the air in the room; $N_{\rm oc}$ is the number of persons occupying the room, each being a source of heat of average power $W_{\rm oc}$ [W]. The system of equations (\ref{eq:phys1}) and (\ref{eq:phys2}) is nonlinear because of the variation of the mass value $m$ over time and the presence of the term responsible for a heat exchange due to ventilation $C_p Q_{\rm in} (T_{\rm in} - T)$. In order to numerically integrate the equations, a first-order explicit method is used to discretize Eqs. (\ref{eq:phys1}) and (\ref{eq:phys2}):
 
\begin{eqnarray}
	\nonumber
	\label{eq:phys1disc}
	m_{i}C_p \frac{T_{i+1}-T_{i}}{\Delta t} & = & U (T_{{\rm out}, i} - T_{i}) +U^{\star} (T^{\star}_{i} - T_{i})+W_{\rm oc}N_{{\rm oc},i}\\
	&+&W_{\rm h/c, j} + C_p Q_{{\rm in},j} (T_{\rm in} - T) + C_p m R_r  (T_{{\rm out},i} - T_{i})\\
	m^{\star}_{i} C^{\star} \frac{T^{\star}_{i+1}-T^{\star}_i}{\Delta t} & = & -U^{\star} (T^{\star}_{i} - T_{i})
	\label{eq:phys2disc}
\end{eqnarray}

\noindent where ${i}$ denotes the current time step, ${i+1}$, the next time step, and ${j}$ the control step. The integration time step ${\bigtriangleup t}$ is chosen to be 1 min so that solution is not changed with its further decrease. The control time step is the interval where controlled parameters, i.e. $(W_{\rm h/c},Q_{\rm in})$ are fixed. For the test cases we present in the next section, the control time step is 1 hour, while the control horizon at which an objective function is minimized is varied: it is set to capture the microclimate dynamics sufficiently accurately while permitting the computation of the optimization problem solution in a reasonable time on a standard computer (less than one minute).

\subsubsection{Weather and occupancy forecasts} 
\begin{figure}[h!]
	\centering
	\includegraphics[width=70 mm]{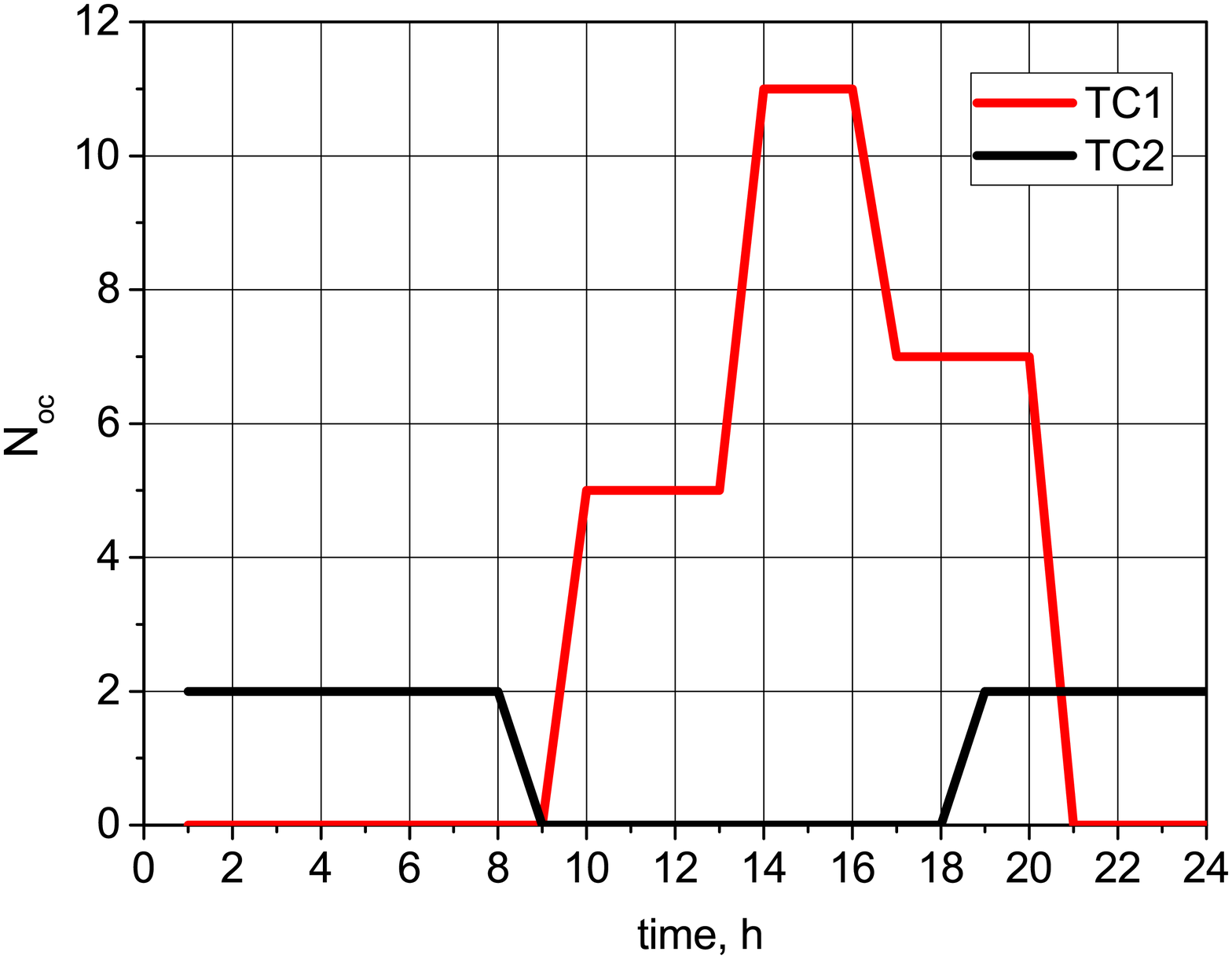}
	\caption{Occupancy variation forecast over 24 hours.}
	\label{fig:Noc}
\end{figure}

\begin{figure}[h!]
	\centering
	\includegraphics[width=70 mm]{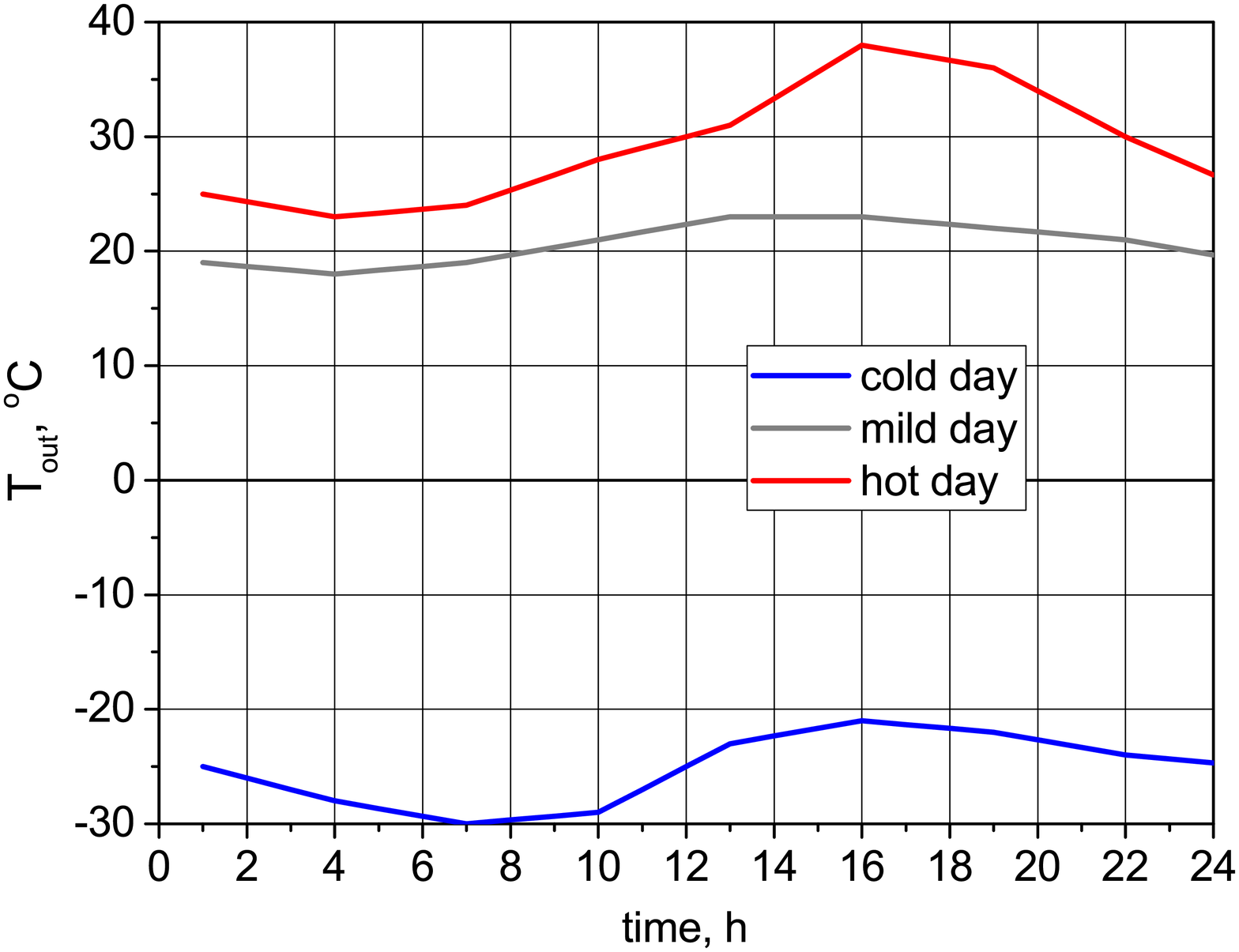}
	\caption{Outside temperature variations forecast over 24 hours corresponding to cold, mild, and hot days.}
	\label{fig:Tout}
\end{figure}

One of the main advantages of MPC is its ability to take into account foreseeable exogenous factors for microclimate dynamics simulations. For the particular cases we address in the present work, these factors are the occupancy of rooms and the external weather. The room occupancy for the test cases is drawn from our knowledge of buildings usage (Fig.~\ref{fig:Noc}): office building occupancy statistics and private house usage based on general understanding of how apartments are used. The weather forecast can, for example, be downloaded from weather forecast websites. Examples of temperature profiles used are shown on Fig.~\ref{fig:Tout} and correspond to generic cold, mild, and hot days in Moscow. Our simulations are thus based on illustrative cases, for which the weather conditions are summarized as either ``typical'' hot, cold or mild days. To define these typical days, we used an on-line weather service and downloaded the weather data for the Moscow region for the last 5 years, with temperatures given every 3 hours for each day recorded; we then calculated averages of the hottest summer days (i.e. days with temperatures above the average over the whole summers), coldest winter days (i.e. days with temperatures below the average over the whole winters) and 100 mild days distributed over the months of May, June, August, and September months each year, during the last 5 years.

\subsection{Constraints} 
Constraints impose the range of HVAC operational regimes and comfort. The lower limit of $W_{\rm h/c}$ represents the air conditioner maximum power (with a negative sign for an energy sink), and the upper limit, the heating maximum power. The ventilation upper limit corresponds to the air handling unit maximum performance. The lower limit depends on the room occupancy: for a fixed number of occupants $N_{\rm oc}$ one need to provide a minimum ventilation rate $Q_{\rm in}^{\rm min}$ in order not to exceed a certain comfort limit of CO$_2$ concentration $\nu_{{\rm CO}_2}^{({\rm max})}$. From the mass balance equation for the CO$_2$ fraction of air, we obtain: 

\begin{equation}
Q_{\rm in}^{\rm min} = \frac{N_{\rm oc} \tilde{Q}_{{\rm CO}_2}}{\nu_{{\rm CO}_2}^{\star} - \nu_{{\rm CO}_2}^{({\rm max})}}
\label{eq:Qmin}
\end{equation}

\noindent where $\tilde{Q}_{{\rm CO}_2}$=0.000012 kg/s is the CO$_2$ mass exhausted by one person each second, $\nu_{{\rm CO}_2}^{\star}$=400 ppm is an average CO$_2$ concentration in the environment and $\nu_{{\rm CO}_2}^{({\rm max})}$ is taken to be 1000~ppm as a maximum comfort value. 

Temperature limits are also important for the comfort: when there are no occupants inside, the lower and upper limits are 15$^{\circ}$C and 25$^{\circ}$C respectively. During daytime the limits are $T^{{\rm comf}}-1 <T<T^{{\rm comf}}+1$ (in $^{\circ}$C). The comfort temperature equals $T^{{\rm comf}}$=22$^{\circ}$C in our cases. Further, to measure a temperature comfort violation, the following penalty function over a 24-hour horizon is introduced in the present work: $T_{\rm penalty}=\sum_{1}^{24}(\delta T\delta t)$ [$^{\circ}$C$\cdot$h], where $\delta t$ is the time when the room temperature takes values outside the comfort interval, with corresponding value $\delta T$.

\subsection{Objective functions}  
In the present work, the objective function $f_{\rm ec}$ represents the energy consumption for a room heating/cooling, ventilation air heating and propulsion:

\begin{equation}
f_{\rm ec} = \sum_{t=1}^{t_{\rm max}} \left[|W_{\rm h/c}(t)| + C_p Q_{\rm in}(t) |T_{\rm in} - T_{\rm out}(t)| + \alpha Q_{\rm in}^3(t)\right]
\end{equation} 

\noindent where $W_{\rm h/c}(t)$ [W] is the heating/cooling power (controlled); $Q_{\rm in}(t)$ [kg$\cdot$s$^{-1}$] is the ventilation mass flow rate (controlled); $T_{\rm in}$ [$^{\circ}$C] is the inflow air temperature; $T_{\rm out}$ [$^{\circ}$C] is the outside air temperature; $\alpha = (2S_{\rm p}\rho)^{-2}$ [m$^2\cdot$kg$^2$] is a proportionality coefficient defined by the product of $S_{\rm p}$ [m$^2$], which is the cross section area of the air inflow device's pipe, and $\rho=1.2$ [kg$\cdot$m$^{-3}$], which is the air density; and $t_{\rm max}$ [h] is the control horizon. The function $f_{\rm ec}$ yields values given in [W$\cdot$h]. It is worth noticing that the function can easily be modified to correspond to a primary energy, total energy or price via multiplication of the function term by a specific constants or functions responsible for, e.g., pricing or efficiency.

\subsection{Optimizer}  
For our purpose here, we use two optimizers, linear and nonlinear, to solve non-convex problems for which global optimatility is difficult if not impossible to prove. 

\subsubsection{Nonlinear non-convex optimization}
The first controller, nonlinear MPC, is based on the nonlinear physical model given by Eqs.~(\ref{eq:phys1disc}) and (\ref{eq:phys2disc}), and requires a nonlinear non-convex optimization algorithm for an optimal control vector calculation. The main advantage of the method is its accuracy provided by the model. We used the optimization solver {\tt fmincon} for nonlinear constrained optimization with nonlinear constraints \cite{Waltz2006}; it is available as part of the Matlab package.

\subsubsection{Linear optimization}
Linear optimization is used to limit the computational burden, so we also developed a linear MPC version of our approach to check its relative accuracy and value as a tool for our work. This naturally entails the linearization of the physical model, Eqs.~(\ref{eq:phys1}) and (\ref{eq:phys2}); to this end, both the mass $m$ of air in the room and the effective temperature $T^{\star}$ are taken as constants at each control time step, assuming that it is sufficient to describe the slow air temperature dynamics, and for the ventilation control, we set $Q_{\rm in}(t)=Q_{\rm in}^{\rm min}(t)$. As a relative change of mass is rather small and the ventilation should work at the minimum possible power to decrease heat exchange with the environment we do not expect a significant decrease of the model accuracy due to the model linearization. An explicit first-order time-stepping is then used with an integration time step equal to the control time step, and we obtain a  linear equality constraints matrix that can be used for optimization computation. We used the {\tt CLP} optimization solver of the OPTI Toolbox package for Matlab. This solver is based on the primal simplex method and allows to solve linear optimization problems. \cite{Forrest2012}

\subsection{Rolling horizon} 
We perform rolling horizon calculations to assess the efficiency of the MPC implemented according to the following sequence of steps: 
\begin{enumerate}
\item Initial values for $T$ and $T^{\star}$ are chosen;
\item The MPC code is run to obtain control vectors at a set control horizon;
\item The first control values are used to simulate $T$ and $T^{\star}$ in 5 to 15 min;
\item The obtained $T$ and $T^{\star}$ are used as initial values for MPC, and the procedure restarts from step 1.
\end{enumerate}

\section{Test cases}
The MPC performance is tested on two cases with the predefined occupancy profiles shown on Fig.~\ref{fig:Noc}. We evaluate the efficiency of MPC by comparison with a standard on/off controller: 
\begin{enumerate}
	\item State space for heating/cooling and ventilation is discretized so that HVAC devices can work at regimes from 0 (switched off) to 10 (maximum power);
	\item Setpoints for temperature correspond to comfort limits used in MPC;
	\item Setpoint for CO$_{2}$ concentration is 900 ppm with a deadband region of 100 ppm;
	\item When no people are inside in hot and normal days the heating/cooling equipment is switched off. 
\end{enumerate}

For each test case we plot the time evolution of the temperature, heating/cooling, and ventilation utilization. 

\begin{table}[h!]
	\centering
	\begin{tabular}{||l c c c||} 
		\hline
		Parameter & Test case 1 & Test case 2 & Test case 2 + supply ventilation \\  
		\hline\hline
		$U$ [W/K]          & 55  & 15 & 15  \\ 
		$U^*$ [W/K]        & 200 & 200& 200 \\
		$m^*C^*$ [MJ]    & 107 & 20 & 20  \\
		$V$ [m$^3$]     & 540 & 105& 105 \\
		$R_{\rm r}$ [h$^{-1}$]    & 0.1 & 0.2& 0.2 \\
		$T_{\rm t=0}$ [$^{\circ}$C]      & 21  & 21 & 21  \\
		$T^*_{\rm t=0}$ [$^{\circ}$C]    & 21  & 21 & 21  \\
		$W_{\rm h/c}^{\rm min}$ [kW] & $-$15 & $-$2 & $-$2  \\
		$W_{\rm h/c}^{\rm max}$ [kW] & 5   & 0.95 & 1.1 \\
		$W_{\rm oc}$ [kW] & 0.12   & 0.12 & 0.12 \\
		$Q_{\rm max}$ [kg$\cdot$s$^{-1}$]   & 0.55& -- & 0.05\\
		$T_{\rm in}$ [$^{\circ}$C]       & 21  & -- & 21  \\
		$S_{\rm p}$ [cm$^2$] & 500  & -- & 120  \\
		\hline
	\end{tabular}
	\caption{Test cases parameters.}
	\label{tab:Params}
\end{table}

\subsection{Test case 1: Skoltech campus room} 
We consider a 100 m$^2$ room which is a standard lecturing room at the Skolkovo Institute of Science and Technology. The relevant physical parameters are derived using the room model assembled in TRNSYS \cite{TRNSYS} and are listed in Table \ref{tab:Params}. The building has a high thermal inertia. When the room occupancy is higher than 20-25 people during lectures, the air quality may deteriorate fairly quickly because of insufficient ventilation. In addition, during a large part of a typical daytime, when no student nor staff is inside, the ventilation system works on its nominal regime thus spending energy, which rather could be saved. 

Rolling horizon simulations are performed for linear and nonlinear MPC. Here, the nonlinear physical model, Eqs. (\ref{eq:phys1}) and (\ref{eq:phys2}) is used to predict the evolution of the microclimate both for the linear (L) and nonlinear MPC. For the LMPC only heating is controlled while the ventilation mass flow ensures that the CO$_{2}$ concentration is less than 1000 ppm. This is confirmed by the nonlinear MPC simulations, where the optimal ventilation mass flow equals the minimum given by Eq. (\ref{eq:Qmin}). 

\subsubsection{MPC vs. LMPC and on/off controller}
For a 24-hour time horizon, the optimization of the vector consisting of 24$\times$3 variables is performed. Here the variables are the temperature, heating/cooling power and ventilation flow for MPC, and the temperature, heating and cooling power only for LMPC to exclude the nonlinearities in the physical model, Eqs. (\ref{eq:phys1disc}) and (\ref{eq:phys2disc}). 

\begin{table}[h!]
	\centering
	\begin{tabular}{||l | c c c | c c c||} 
	  \hline
		&\multicolumn{3}{c}{Heating($+$)/Cooling ($-$) [kWh]} \vline & \multicolumn{3}{c}{Ventilation [kWh]} \vline \\
		\hline
		              Day type & MPC & LMPC  & on/off & MPC & LMPC  & on/off \\ 
	    \hline\hline
		Cold       & 36.8  & 34.6   & 62.0  & 61.7 & 61.7  & 70.4 \\ 
		Mild     & $-$2.8  & $-$3.9   & $-$6.1  & 2.0  & 2.0   & 2.4  \\
		Hot        & $-$11.4 & $-$13.9  & $-$16.1 & 18.7 & 18.7  & 21.7 \\
		\hline
	\end{tabular}
	\caption{Energy consumption for 24 hours of operation for the test case 1: MPC vs. LMPC vs. on/off control.}
	\label{tab:TC1_WQ}
\end{table}

\begin{table}[h!]
	\centering
	\begin{tabular}{||l | c c c | c c c||} 
		\hline
		&\multicolumn{3}{c}{$T_{\rm penalty}$ [K$\cdot$h]} \vline  \\
		\hline
		Day type & MPC & LMPC  & on/off \\ 
		\hline\hline
		Cold       & 0.0  & 13.1   & 0.2  \\ 
		Mild     & 0.0  & 0.0    & 0.0  \\
		Hot        & 0.0  & 0.1    & 0.2  \\
		\hline
	\end{tabular}
	\caption{Temperature comfort violation for TC1.}
	\label{tab:TC1_Penalty}
\end{table}

\begin{figure}
	\centering
	\includegraphics[width=122 mm]{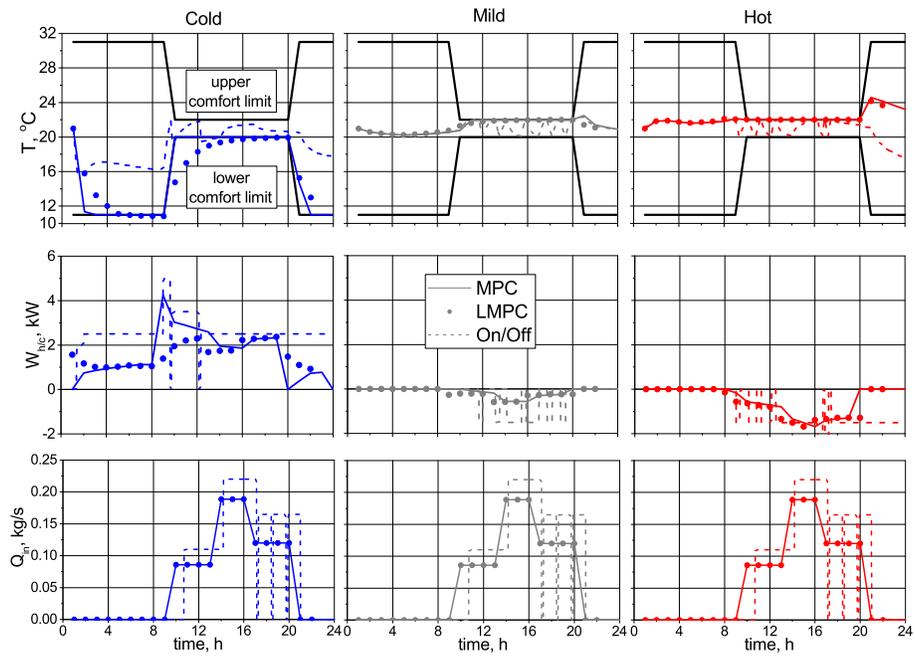}
	\caption{TC1. Comparison of MPC (solid line), LMPC (dots) and On/Off control (dashed line) for different weather conditions: cold (left column), mild (central column) and hot (right column).}
	\label{fig:Cold_MPCvsLMPC}
\end{figure}

On Fig. \ref{fig:Cold_MPCvsLMPC}, one can see from the temperature time evolution, that MPC allows a proper comfort level for harsh conditions while the LMPC underestimates the necessary heating level. The corresponding values of $T_{\rm penalty}$ are given on Table~\ref{tab:TC1_Penalty}. For both control methods the ventilation flow equals the lower limit. From the results displayed in Table~\ref{tab:TC1_WQ}, one can conclude that MPC is the most efficient in terms of energy spent for microclimate with the on/off controller being the worst among the tested approaches.

\subsubsection{Forecast horizon influence} 
\begin{figure}
	\centering
	\includegraphics[width=130 mm]{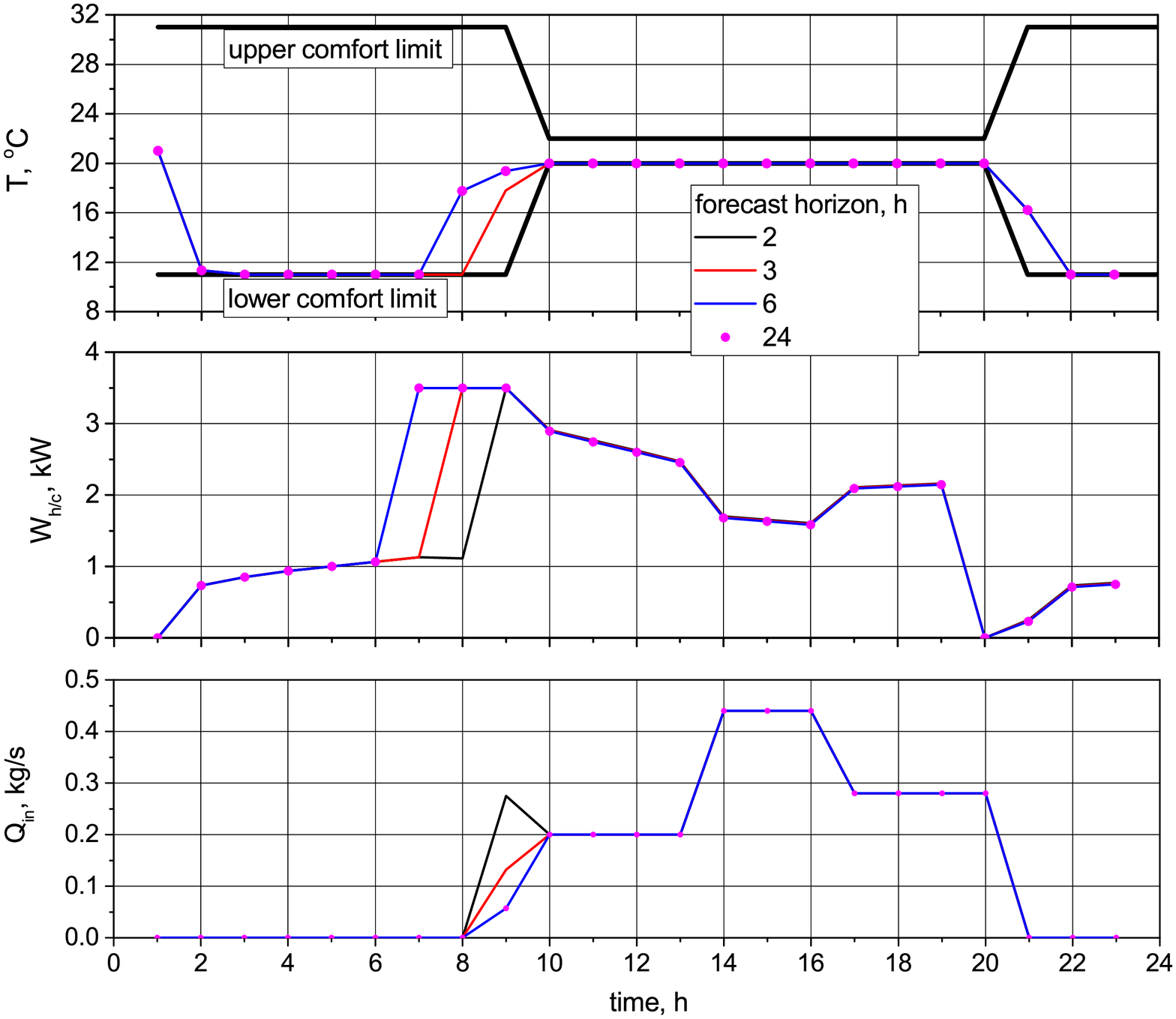}
	\caption{MPC control with different prediction horizons under cold weather. Maximum heating power = 3.5 kW.}
	\label{fig:Cold_2h_3h_4h_6h}
\end{figure}

To obtain the optimal time horizon for MPC, simulations of MPC with different horizons are performed. In practice the power of heating/cooling systems is limited and hence one cannot rapidly heat/cool the room to a desired level. Therefore heating/cooling equipment should be switched in advance. For example, the numerical results depicted on Fig. \ref{fig:Cold_2h_3h_4h_6h} are obtained for a cold weather with maximum heating power limited to 3.5 kW. Here, the control corresponding to the minimum of the objective function starts from horizon = 4 hours and no further change of control is observed. The corresponding objective functions are 184 kWh for horizon = 2 hours, 181 kWh for 3 hours and 180 kWh for larger horizon values.

\subsubsection{Uncertainty influence}
In practice, as occupancy and weather forecasts are not perfectly accurate, errors in future states estimation are unavoidable. To account for this problem at every rolling horizon step, stochastic errors are added to initial values of occupancy, outside temperature, and room temperature. The disturbances are introduced as follows with means $\mu$ and standard deviation $\sigma$:

\noindent$T^{\rm (initial)}$: added normally distributed noise with $\mu=0$ and $\sigma=1$ K;\\
$N_{\rm oc}^{\rm (initial)}$: multiplied by uniformly distributed value from 0 to 2;\\
$T_{\rm out}^{\rm (initial)}$: added normally distributed noise with $\mu=0$ and $\sigma=1$ K.

\begin{figure}[h!]
	\includegraphics[width=130 mm]{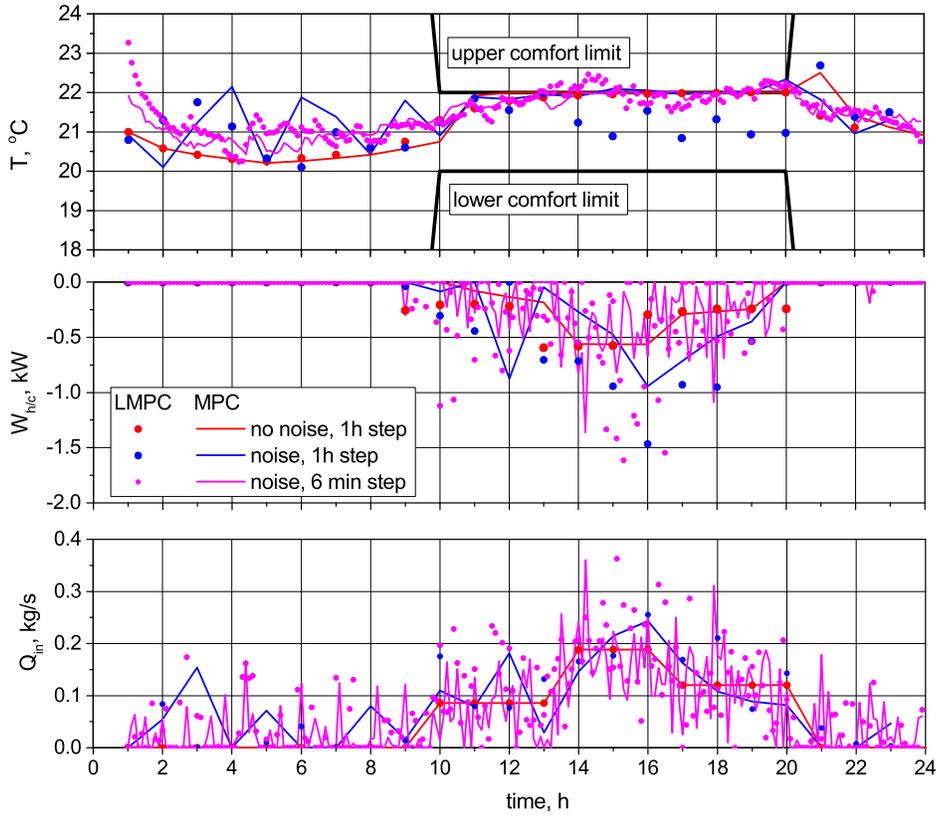}
	\caption{Time evolution of temperature, heating/cooling power, and ventilation mass flow rate accounting for uncertainty, for the mild weather case. MPC (lines) and LMPC (dots).}
	\label{fig:Norm_T_dist}
\end{figure}

\noindent Simulations are performed with MPC and LMPC for a normal day (see Fig. \ref{fig:Tout}). Different discretization time steps in rolling horizon are used: 1 hour, and 6 minutes. From the results displayed on Fig. \ref{fig:Norm_T_dist}, one can see that in general both MPC and LMPC provide rather good control in terms of comfort. But LMPC violates comfort boundaries more than MPC due to the loss accuracy of the physical model used. Corresponding temperature penalties values equal $T_{\rm penalty}^{\rm (MPC)} \approx 0.1$ and $T_{\rm penalty}^{\rm (LMPC)} \approx 0.8$.

\begin{figure}
	\includegraphics[width=130 mm]{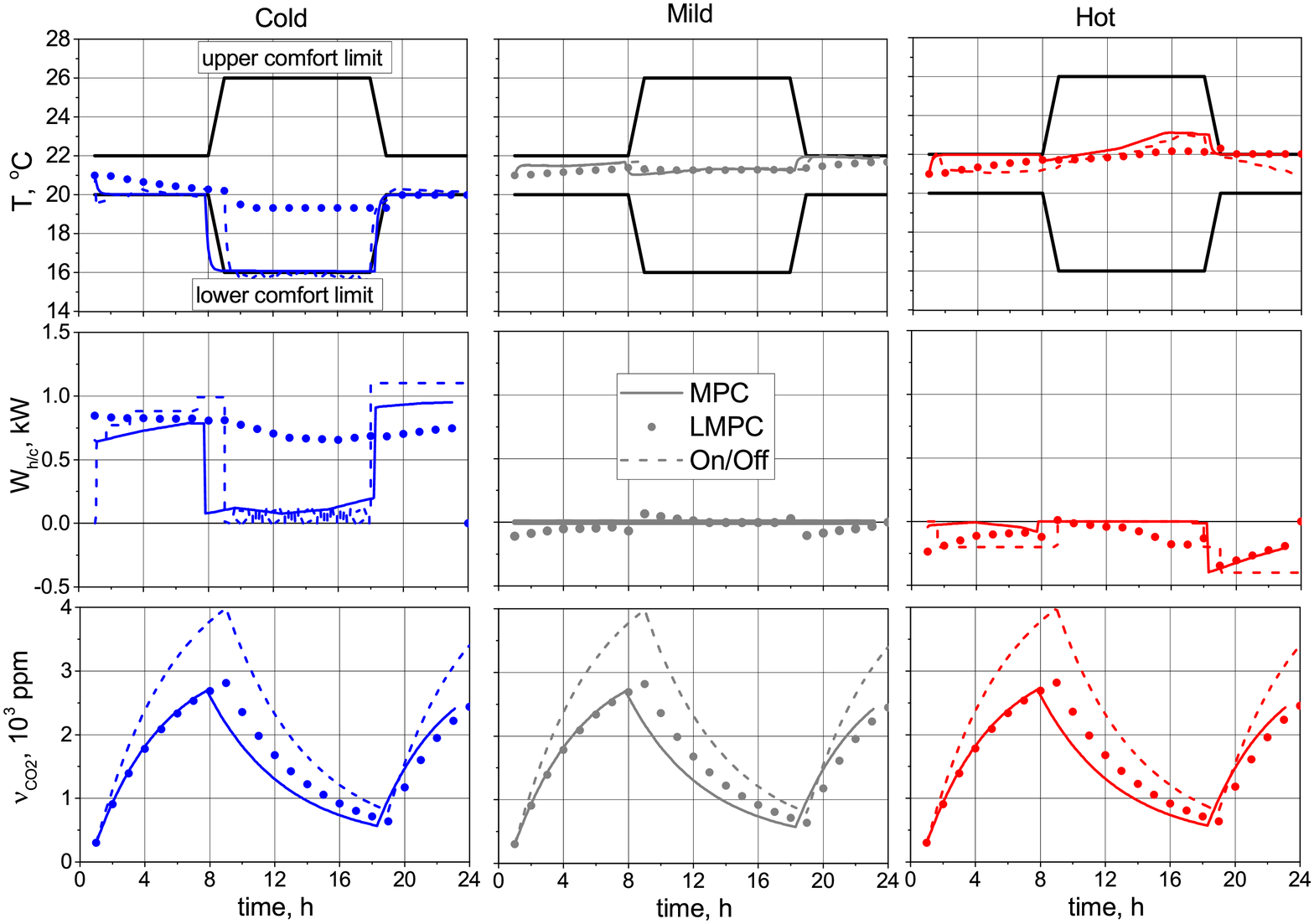}
	\caption{TC2 without a supply ventilation system. Comparison of MPC (solid line), LMPC (dots) and On/Off control (dashed line) for different weather conditions: cold (left column), mild (central column) and hot (right column).}
	\label{fig:TC2nobr}
\end{figure}

\begin{figure}
	\includegraphics[width=130 mm]{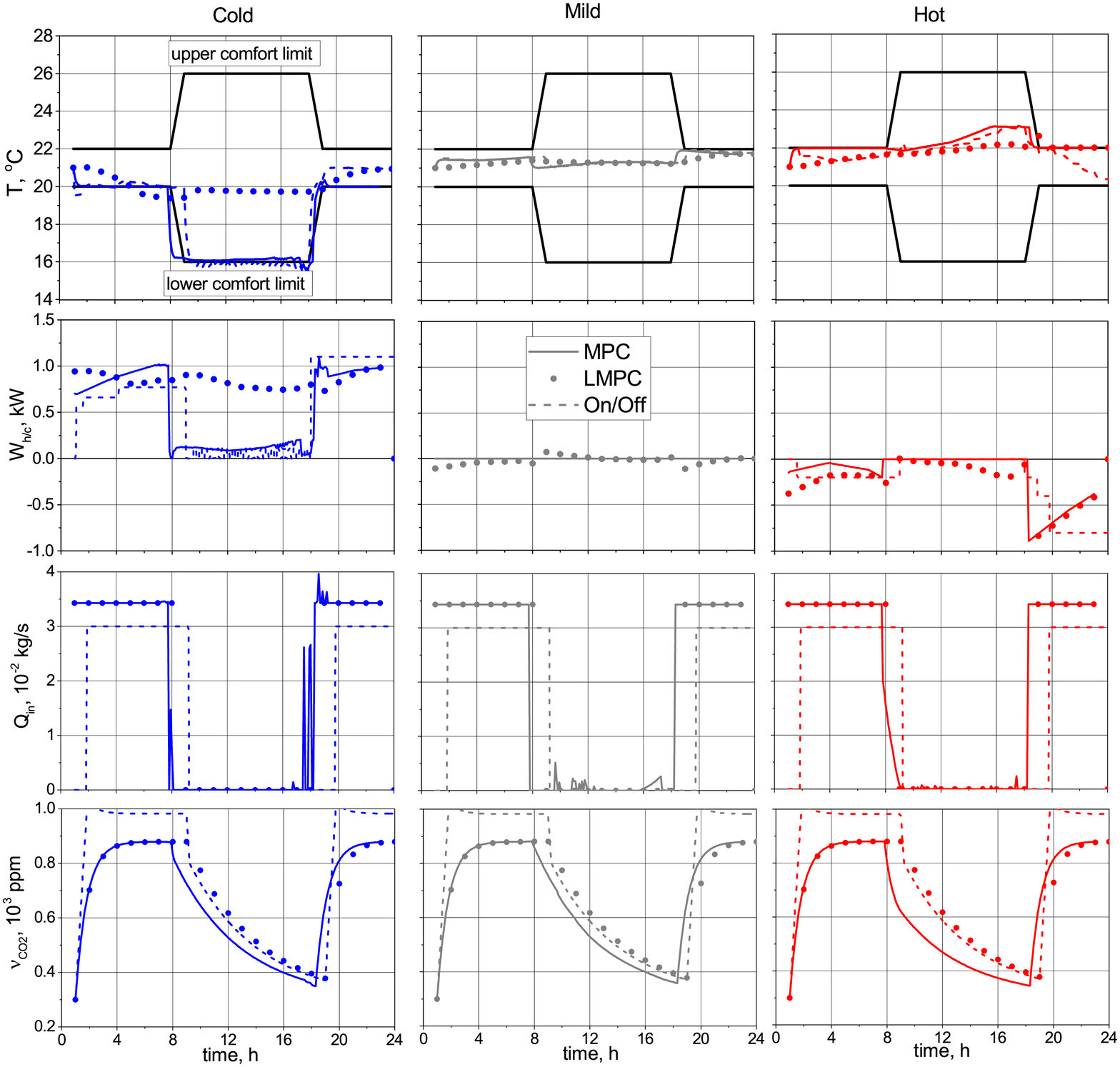}
	\caption{TC2 with a supply ventilation system. Comparison of MPC (solid line), LMPC (dots) and On/Off control (dashed line) for different weather conditions: cold (left column), mild (central column) and hot (right column).}
	\label{fig:TC2br}
\end{figure}

\subsection{Test case 2: Single-room apartments} 
The second test case corresponds to a typical single-room apartment (Table \ref{tab:Params}). Buildings series P44 were the most popular panel house type build in 1980-1990, and up to date they are among the most widespread in some cities. One of the main drawbacks of these and most of other apartments is the absence of a proper ventilation system which together with low infiltration leads to inappropriate growth of carbon dioxide concentration when occupants are inside for a long period of time, e.g. at night (see Fig. \ref{fig:TC2nobr}). Formerly these apartments were equipped with rather low-quality windows not well fitted, thus leading to high air leakage; but in recent decades these outdated windows started to be replaced with modern leakproof plastic glass units that significantly decreased the overall infiltration rate of apartments. To allow fresh air to get into a room one can open a window. But in megacities like Moscow this also causes noise, excessive heat leakage, as well as pollutants like NO$_x$ and particulate matter penetration. This solution thus generates uncontrolled air flows in a flat, which may yield discomfort, nay health problems. 

\begin{table}[h!]
	\centering
	\begin{tabular}{||l |c c c |c c c||} 
	  \hline
		&\multicolumn{3}{l}{Heating(+)/Cooling (-), [kWh]} \vline & \multicolumn{3}{l}{Ventilation, [kWh]} \\
		\hline
	 	Configuration                & MPC   & LMPC  & on/off & MPC   & LMPC & on/off \\ 
		\hline\hline
		Cold day, with SVS       & 12.2  & 19.5   & 13.8  & 18.8   & 19.7  & 18.1  \\
		Cold day, no SVS         & 11.1  & 17.2   & 12.9  & --     & --    & --    \\ 
		Mild day, with SVS     & 0     & 0      & 0     & 0.6    & 0.6   & 0.6   \\
		Mild day, no SVS       & 0     & 0      & 0     & --     & --    & --    \\
		Hot day, with SVS        & $-$3.8  & $-$5.7   & $-$5.7  & $<$0.1 & $<$0.1& $<$0.1\\
		Hot day, no SVS          & $-$1.7  & $-$3.2   & $-$4.1  & --     & --    & --    \\
		\hline
	\end{tabular}
	\caption{Energy consumption for 24 hours of operation for the test case 2: MPC vs. on/off control.}
	\label{tab:TC2_WQ}
\end{table}

\begin{table}[h!]
	\centering
	\begin{tabular}{||l | c c c | c c c||} 
		\hline
		&\multicolumn{3}{c}{$T_{\rm penalty}$ [K$\cdot$h]} \vline     \\
		\hline
		Day type                   & MPC  & LMPC   & on/off \\ 
		\hline\hline
		Cold day, with SVS     & 0.4  & 1.7    & 0.2    \\ 
		Cold day, no SVS       & 0.4  & 0.7    & 0.1    \\ 
		Mild day, with SVS   & 0.0  & 0.0    & 0.0    \\
		Mild day, no SVS     & 0.0  & 0.0    & 0.0    \\
	    Hot day, with SVS      & 0.0  & 0.7    & 0.0    \\
	    Hot day, no SVS        & 0.0  & 0.4    & 0.0    \\
		\hline
	\end{tabular}
	\caption{Temperature comfort violation for TC2.}
	\label{tab:TC2_Penalty}
\end{table}

Window opening is neither an energy efficient nor a comfort solution to deal with fuggy indoor environments. A centralized inflow ventilation system installation may be considered to the appropriate solution, but normally it is not for most of the popular series dwellings because of limited ceiling height. Another solution for house is the supply ventilation system, which is increasingly popular in Russia: a small box installed over an orifice on the wall that includes a compressor, filters and a heater. The maximum power of these devices is typically 1.5 kW, most of which being spent for the inflow air heating. SVS are also equipped with a filtering (NO$_x$, SO$_x$, etc.) and disinfecting system (viruses, bacteria). On Fig. \ref{fig:TC2br} one can see that the solution requires higher heating/cooling power available but this allows to limit carbon dioxide concentration at a desired level. Its advantages over windows opening are controlled fresh air mass flow rate that is needed for particular indoor conditions, air cleaning and outside noise reduction. As seen from table~\ref{tab:TC2_WQ} in the present test case, the MPC does not seem to be a necessary control solution as it neither allows to save energy nor is it able to improve the comfort level. The reason can be the quite limited thermal inertia of the room tested in TC2. From Table~\ref{tab:TC2_Penalty} one can conclude that all the control methods provide almost the same comfort level. However, if energy storage is going to be utilized or a significant change in prices for energy (both electric and thermal) are to be introduced MPC can become the solution able to significantly improve a microclimate efficiency.

\section{Concluding remarks} 
The scope of our work is timely: as the majority of dwellings are highly dependent on fossil fuels for the provision of energy, either a thorough renewal of the housing stock or even their complete renovation up to nowadays best standards are prohibitive in terms of costs. Further, in a context where renewable energy sources cannot yet be seen as full substitutes for fossil fuels, even the strategies which promote low-power consumption and environment-friendly technologies (or ``best available technologies'' - BAT) have a major drawback known as the rebound effect \cite{FontVivanco2016}: as these BAT multiply to satisfy the ever growing public demand, the toll on standard energy sources is becoming unsustainable, especially when indoor comfort acquires a central role in people's daily life. Therefore on a short to mid-term time scale, only smart technologies can answer the challenges of quality indoor microclimate at a reasonable cost through home automation systems and control. 

In this article, we studied how to control comfort and efficiency of energy utilization in buildings. The two test cases we presented serve different purposes, one is modern, one dates back to the early 1980s; their thermal and ventilation characteristics are quite different but both are subjected to the same high-amplitude yearly variations of the weather: twice a year the outside temperature undergoes a 70$^{\circ}$C change. This has a strong impact on the buildings, but also on the indoor microclimate needs, which require proper strategy definition for retrofitting where needed. Further, while the two actual test cases we treat are located in Moscow, our results and approach naturally extend to places, which have a similar climate. To contribute a solution, we developed a method based on model predictive control. MPC simulations for the modern building case study show the significant advantage that nonlinear MPC has over both linear MPC and standard on/off controller in providing comfort with lower energy consumption, coupled with the ability to operate the HVAC equipment in advance in order to prepare the desired microclimate at a given time; in other words MPC can utilize the capacity for shaving peaks in demand via pre-cooling or pre-heating the room at off-peak hours. In particular, for hot and cold seasons MPC allows to provide significantly higher levels of comfort than LMPC solutions. In addition, solutions which provide a good indoor environment also mitigate the detrimental effects of harsh weather conditions.

With the second test case corresponding to an old building, we showed that this type of apartments requires a ventilation system to be installed as, under standard conditions of use, the air quality is rather poor with, in particular, CO$_2$ concentrations much higher than what is advised by comfort norms. MPC for the HVAC control here, does not appear to provide acceptable solutions as we observe no significant reduction of energy consumption. This implies that direct application of one type of solutions, successful in some cases, is not a valuable strategy as the specific characteristics of buildings must be assessed and solutions adapated; this conclusion is in line with recent work putting forth the need for suitability assessment \cite{Geng2015}. Further, contextual factors influence the investment strategy for building or house management \cite{Nair2010}. Nevertheless, MPC can be a promising solution if a storage device is utilized as part of the whole system. The development of a multi-zone MPC controller with an adaptive data-driven approach to the microclimate dynamics model description and its practical implementation constitutes the next step of the present work.

\section*{Acknowledgments}
This work is supported by the Skoltech NGP Program (Skoltech-MIT joint project). 

\section*{References}

\bibliographystyle{elsarticle-num}

\end{document}